\documentclass[preprint,12pt]{elsarticle}

\usepackage{amssymb}

\usepackage{dsfont}

\usepackage{eurosym}

\usepackage[utf8]{inputenc} 
\usepackage[T1]{fontenc}    
\usepackage{hyperref}       
\usepackage{url}            
\usepackage{booktabs}       
\usepackage{nicefrac}       
\usepackage{microtype}      
\usepackage{lipsum}		
\usepackage{graphicx}
\usepackage{natbib}
\usepackage{doi,amssymb,amsthm}
\usepackage{amsmath}

\usepackage{comment}

\usepackage{algorithm}
\usepackage{algpseudocode}

\journal{Energy Economics}

\begin{document}

\begin{frontmatter}



\title{Reinforcement Learning for Bidding Strategy Optimization in Day-Ahead Energy Market}


\author[inst1]{Luca Di Persio}

\affiliation[inst1]{organization={University of Verona, Department of Computer Science},
            addressline={Strada le Grazie 15}, 
            city={Verona},
            postcode={37134}, 
            country={Italy}}

\author[inst1]{Matteo Garbelli}
\author[inst2]{Luca Maria Giordano}

\affiliation[inst2]{organization={University of Milano, Department of Mathematics},
            addressline={Via Cesare Saldini 50}, 
            city={Milano},
            postcode={20133}, 
            country={Italy}}

\begin{abstract}
In a day-ahead market, energy buyers and sellers submit their bids for a particular future time, including the amount of energy they wish to buy or sell and the price they are prepared to pay or receive. However, the dynamic for forming the Market Clearing Price (MCP) dictated by the bidding mechanism is frequently overlooked in the literature on energy market modelling. Forecasting models usually focus on predicting the MCP rather than trying to build the optimal supply and demand curves for a given price scenario. Following this approach, the article focuses on developing a bidding strategy for a seller in a continuous action space through a single agent Reinforcement Learning algorithm, specifically the Deep Deterministic Policy Gradient. The algorithm controls the offering curve (action) based on past data (state) to optimize future payoffs (rewards). The participant can access historical data on production costs, capacity, and prices for various sources, including renewable and fossil fuels. The participant gains the ability to operate in the market with greater efficiency over time to maximize individual payout.
\end{abstract}


\begin{highlights}
\item We propose a Reinforcement Learning algorithm to optimize bidding strategies for electricity suppliers.

\item The supplier interacts with the market model by submitting offering curves based on historical electricity prices. The feedback received as a reward is used to refine the bidding strategy.

\item We calibrate our model with real-world data. Since the action space is high-dimensional and continuous, we employ the Deep Deterministic Policy Gradient algorithm to approximate the optimal policy and the value function.

\end{highlights}


\end{frontmatter}



\section{Introduction}
\label{sec:intro}

The price of electricity in the European market can fluctuate significantly due to the varying modes of production from different sources, each with its constraints related to weather, production volume, and storage challenges. From a financial, economic, and ecological point of view, modelling market movements is crucial to maximizing power output and helping the energy transition process.

\medskip

Throughout the article, we focus on a specific class
of electricity trading markets: day-ahead energy markets. In general, these markets use a simple auction paradigm: operators on the supply and demand side submit bids that include a quantity $(q)$ and a price per unit $(p)$ for the trading period(s) (typically hours) of the following day. In this basic setup, the intersection point of the demand and supply curves is sought, which ensures a balance of consumption and production while determining the Marginal Clearing Price (MCP). A market clearing algorithm called EUPHEMIA (Pan-European Hybrid Electricity Market Integration Algorithm) \cite{euphemia} has been developed among European power exchanges to facilitate the integration of European power markets.

\medskip

The article focuses on developing an algorithm to solve the stochastic optimal control problem related to a market participant operating in a day-ahead power market whose equilibrium, for most European countries, is computed via the Euphemia algorithm. The optimization problem is solved via Reinforcement Learning (RL), where an agent, i.e. an energy operator, interacts with a stochastic environment whose state is described by the historical values of electrical price to maximize its profit. Rather than directly looking at its profit, the agent selects a proper offering curve and receives a reward based on the new state of the environment, i.e. the prices for the next day used to compute the reward. Our purpose concerns computing the deterministic, offline optimal policy that produces the best action to maximize the cumulative discounted reward given a state of historical prices.

\medskip

The paper is organized as follows: we conclude Sec. \ref{sec:intro} by a review of RL methodology and its application to energy market problems; 
 in Sec. \ref{sec:model} we introduce the theoretical setting for stating the electricity auction problem as an optimal control one; in Sec. \ref{sec:ddpg} we present the Deep Deterministic Policy Gradient (DDPG) method and its adaptation to the electricity auction framework; in Sec. \ref{sec:numerics} we present the results of the numerical simulations and report some considerations and limitations of the algorithm. We conclude the article with Sec. \ref{sec:future_directions} sketching future directions that may employ this algorithm as a reference starting point.

 \medskip

{\em Literature review}. RL \cite{SB,seb} is a learning paradigm that maps situations to actions to maximize a numerical reward signal through repeated experience gained by interacting with the environment. The agent aims to develop a strategy that maximizes the expected cumulative reward over time by learning a policy that maps states to actions. Some of the most common algorithms for RL rely on learning optimal action-value functions by computing the corresponding Q-value, i.e. the quality, the optimal expected future value of the selected action given a particular space. Reference works for Q-learning are contained in \cite{watkins1992q}. For its extension, the Deep Q-Networks (DQN), we refer to \cite{mnih2015human}. A further step is introduced by Actor-criticism methods, which combine the advantages of policy gradient methods and value function approximation to improve the learning process. The actor is responsible for generating actions based on the current policy, while the critic learns to evaluate the policy by estimating the value function. We refer to \cite{actorcritic} for a complete discussion.

DDPG is an off-policy algorithm that extends the idea of the actor-critic method to continuous action spaces \cite{lil}. DDPG uses a deep neural network to approximate the policy and another deep Neural Network (NN) to approximate the value function. Throughout our research, we use the RL model developed in \cite{lil} as the reference algorithm for our setting by developing its adaptation for the setting considered. Following this track, DDPG has been applied to learn optimal bidding strategies for generators and energy storage systems in day-ahead markets and real-time markets \cite{ddpg2}. Focusing on other projects that applied RL algorithms to the energy field, we start by citing \cite{xho}: the authors model the electricity auction market using a $Q$-learning algorithm considering each supplier bidding strategy as a Markov Decision Problem where the agents learn from experience an optimal bidding strategy to maximize its payoff. Although there are certain limits in terms of application - in the case studies considered in \cite{xho} such as the use of simple synthetic datasets as well as discrete Q-tables for the pairings of actions-state, this work serves as a landmark for the research developed in this article. The main differences with the problem studied in \cite{xho} rely on the source of data since we choose to employ historical times series rather than using synthetic data such as in \cite{xho} as well as the development of a more complex RL model. Q-learning has always been used in electricity auctions to learn bidding strategies for market participants, such as generators and retailers \cite{QlearningAuction}. However, the discrete nature of Q-learning can limit its applicability to auctions with large or continuous state and action spaces. In contrast, policy gradient methods can handle continuous state and action spaces, making them suitable for electricity auctions with complex market dynamics. One limitation of policy gradient methods is that they may require a large batch of samples for stable learning. 
Outside the applications to the day-ahead market, we can find, e.g. \cite{actor}, where actor-critic methods are applied for learning bidding strategies and demand response for regulating the power-exchange bidding mechanism in deregulated power system management, offering a balance between exploration and exploitation.

\medskip

\section{Problem Statement and Model Formulation}
\label{sec:model}

\begin{figure}[H]
     \centering
     \includegraphics[scale=0.28]{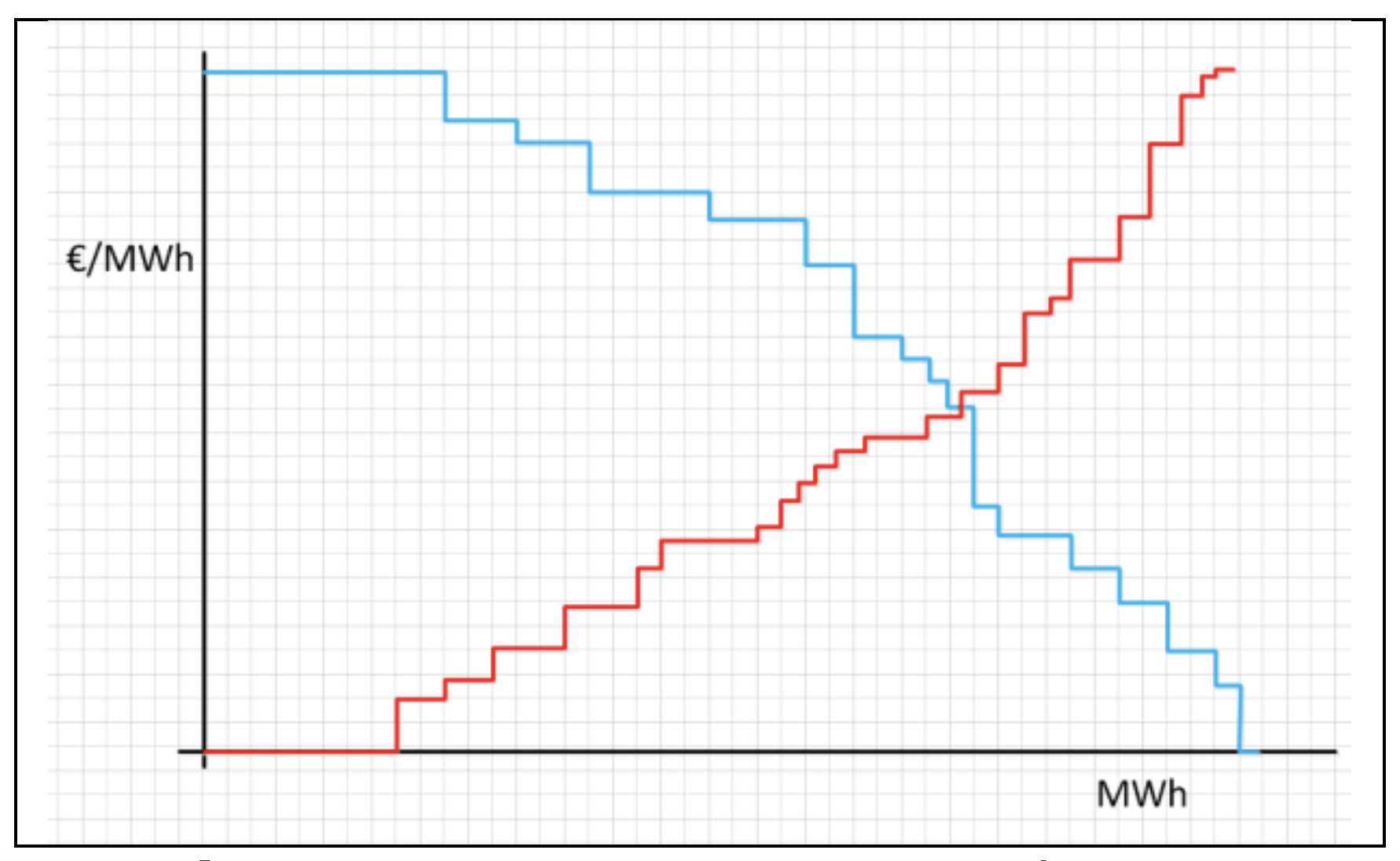}
     \caption{Demand and Offer stepwise aggregated curve (from  \cite{euphemia}).}
     \label{fig:agg_step}
 \end{figure}

Modern energy networks rely heavily on day-ahead electricity markets, which effectively offer a framework to purchase electricity for next-day delivery. Throughout the paper, we focus our attention on the Italian electricity market. The unitary MCP of energy, also known as PUN ({\it Prezzo Unitario Nazionale}) in the Italian market, is established by a bidding procedure involving several market participants, including producers, consumers, and traders. PUN  represents the cost of producing the last unit of electricity needed to meet demand. The Euphemia algorithm sets it \cite{euphemia} as the intersection of the supply and demand curves (see Fig. \ref{fig:agg_step}).

A critical issue for electricity suppliers is how to optimally bid on the auction market to maximize their profit. We model the optimal control problem into a single-agent RL setting, solved via the DDPG algorithm, introduced in \cite{lil}, that uses deep NNs to approximate policy and value functions in a high-dimensional, continuous action space.

\subsection{Stochastic Optimal Control for Electricity Auction Problem}


The producer operates in a stochastic environment for selecting the best bidding strategy, striving to maximize their profit over the long term while meeting the needs of the available resources. At each stage $t$, given a state $s_t \in S$, the seller selects and executes an $a_t \in A$ that, in our setting, corresponds to a stepwise energy/price function, depending on the learned policy $\pi: S \rightarrow P(A)$. The agent's goal is to take actions that will maximize its expected long-term performance with an unknown transition function $P$. To achieve this, the agent learns a behaviour policy $\pi: S \rightarrow P(A)$ that optimizes its expected performance in the long run. The system progresses from state $s_t$ under joint action $a_t \in A$, based on the transition probability function $P$, to the next state $s_{t+1}$, providing updated information on the aggregated load, corresponding unit loads and the new unit price.

\begin{figure}[H]
    \centering
    \includegraphics[width=.69\textwidth]{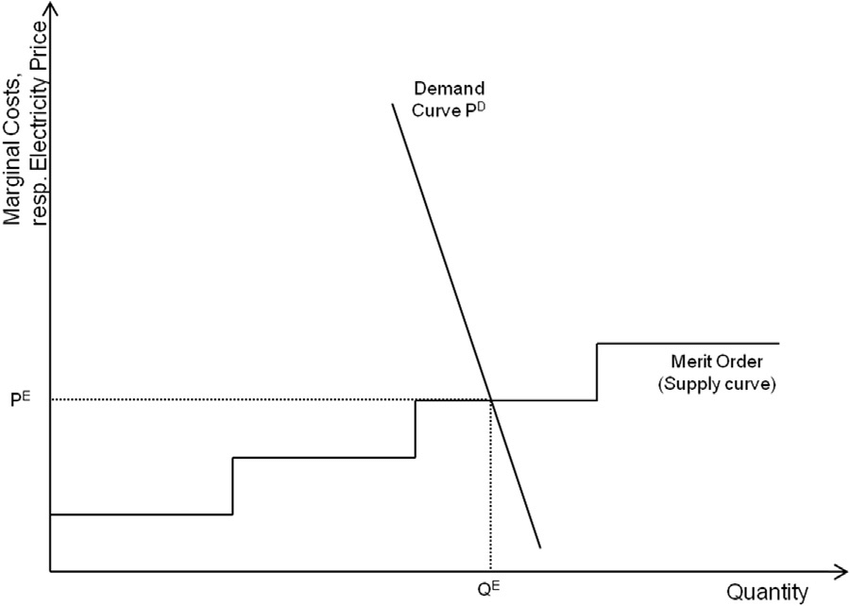}
    \caption{Offering Curve: the intersection with the demand (thus the reward) is computed with the next-day price of electricity.}.
    \label{fig:step}
\end{figure}

In Fig. \ref{fig:step}, we report an example of an offering curve that we obtain corresponding to an action that the agent can take at time $t$.

The reward function $ r_t : S \times A \times S \rightarrow \mathbb{R}$ corresponds to the received feedback signal transitioning from $(s_t, a_t)$ to $s_{t+1}$.
Since the immediate reward is insufficient for providing insights into the long-term profit, it is crucial to introduce the return value $R_t$, defined over a finite time horizon $T$. The following expression gives the return value $R_t$:

\begin{equation}
\label{discount}
\begin{array}{c}
\displaystyle R_t= r_{t+1}(s_t, a_t) + \gamma^{t+1} r_{t+2} (s_{t+1}, a_{t+1}) + \, ... \, + \gamma^{T-1} r_T (s_{T-1}, a_{T-1})= \medskip \\ \displaystyle \qquad =  r_{t+1}(s_t, a_t) + \sum_{i=t+1}^{T-1} \gamma^{i} r_{i+1}(s_i, a_i) \medskip \, .
\end{array}
\end{equation}

In Eq. \eqref{discount}, $\gamma^i \in [0,1]$ corresponds to the discount factors that determine the importance of future rewards compared to the immediate ones, with lower values focusing on short-term rewards. Hence, the return value $R_t$ corresponds to the discounted sum of future rewards, allowing agents to optimize their actions for long-term profit. 
Each agent receives $R_{t+1}$ as immediate feedback for the state transition. Hence, the agent aims to optimize an objective corresponding to the return value $R_t$ 

\begin{equation}
\label{functinoal}
J = \mathbb{E}_{r_t, s_t, a_t \sim \pi} [R_t] \,  ,
\end{equation}

That corresponds to learning a policy that maximizes the cumulative future payoff to be received starting from any given time $t$ until the terminal time $T$.

The agent’s value function associated with such a control problem reads.

\begin{equation}
V ( s_t) = \max_{a_t \in A} \mathbb{E} \left[ R (s_t, a_t, s_{t+1})  \right] \, . 
\end{equation}

The dynamic programming principle implies that $V$ satisfies the Bellman equation.

\begin{equation}
V(s_t) = \max_{a_t \in A} \mathbb{E} \left[  r (s_t, a_t, s_{t+1}) + \gamma V (s_{t+1})\right] \, . 
\end{equation}

In RL, it is useful to define an action-value function to measure the 'quality' of taking a specific action $a$,
 and hence is called the Q function $Q: S \times A \rightarrow \mathbb{R}$ defined as

\begin{equation}
\label{q-values}
Q(s_t,a_t) = R (s_t, a_t) + \gamma \max_{a \in A} Q (s_{t+1}, a_t) \, .
\end{equation}

The $Q$-function solves the Bellman equation by recursively updating the value of taking a specific action 
$a_t$ in a given state $_t$

\begin{equation}
    \label{bellman}
Q(s_t,a_t) = R (s_t, a_t) + \gamma \max_{a_{t+1}} Q (s_{t+1}, a_{t+1}) \, . 
\end{equation}

 By iteratively applying this update, the $Q$-function converges to the optimal action-value function, effectively solving the Bellman equation. This process allows the algorithm to learn the optimal bidding strategy by estimating the quality of different offering curves based on historical price data. We build the critic network to approximate the equation in Sec. \ref{sec:ddpg}. \eqref{bellman}.

\subsection{Construction of the model}
\label{sec:model2}

Following the paradigm of a repeated day-ahead electricity auction market, the agent will attempt to maximize its profit in the long run in a recursive way.

\medskip

At each time step $t$, the agent
receives an observation consisting in of $24$-hour PUN array prices $PUN_t$ of the last $d$ days. We assume a fully observable environment where the state of the environment is represented by the market electricity prices expressed in €$/ MWh$.

\begin{figure}[H]
    \centering
    \includegraphics[width=.55\textwidth]{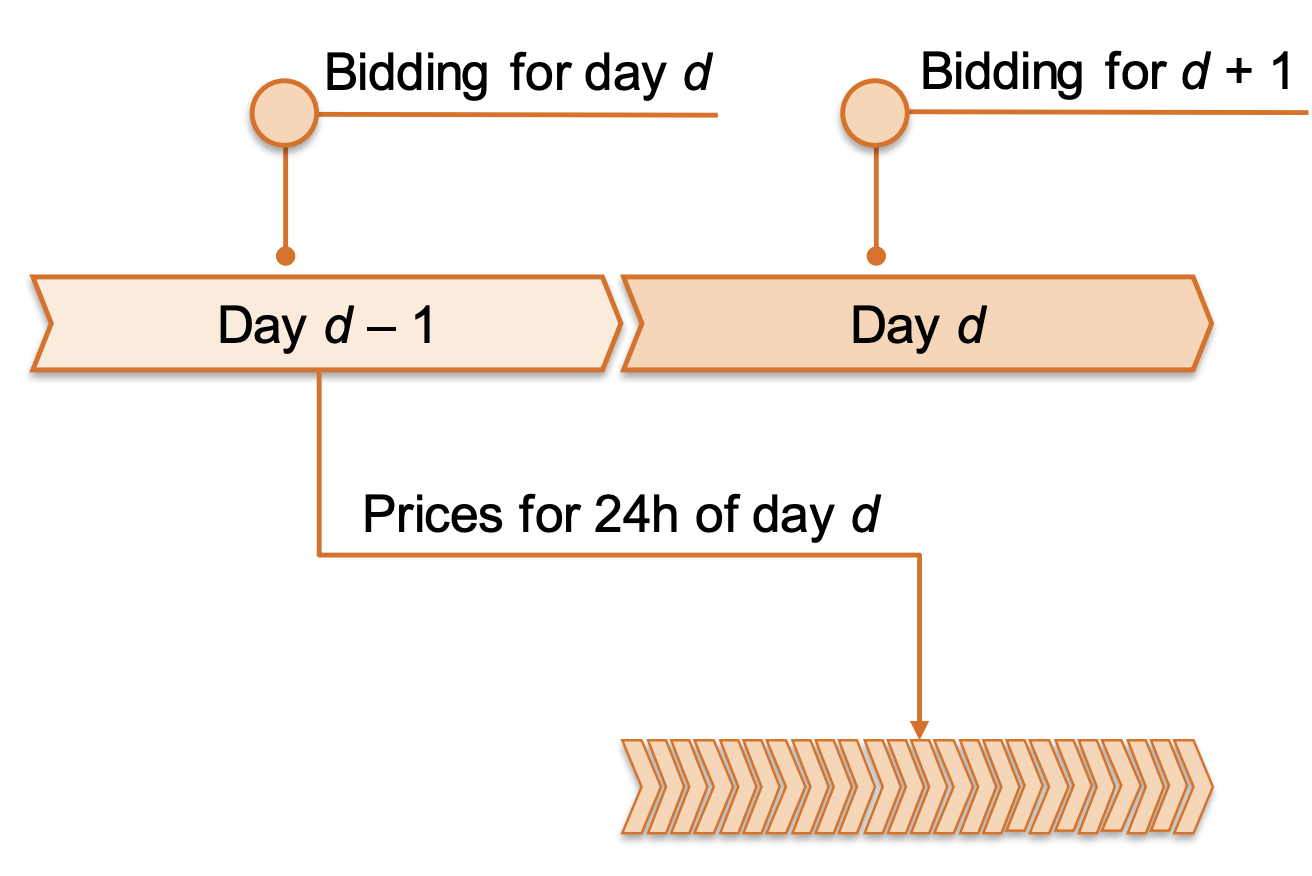}
    \caption{Bidding settlement in day-ahead auctions (from \cite{muw})}.
    \label{fig:bidding}
\end{figure}

We briefly describe how state, actions, and reward are defined for this problem.

\textbf{State}. The current state $s_t$ of the system is given by:

    - \textit{Electricity prices}. An extracted 7-day batch of historical electricity prices is denoted as $PUN_t$, represented by a matrix of $168$ values, $24$ hours times $7$ days.

\begin{figure}[H]
    \centering
    \includegraphics[width=.89\textwidth]{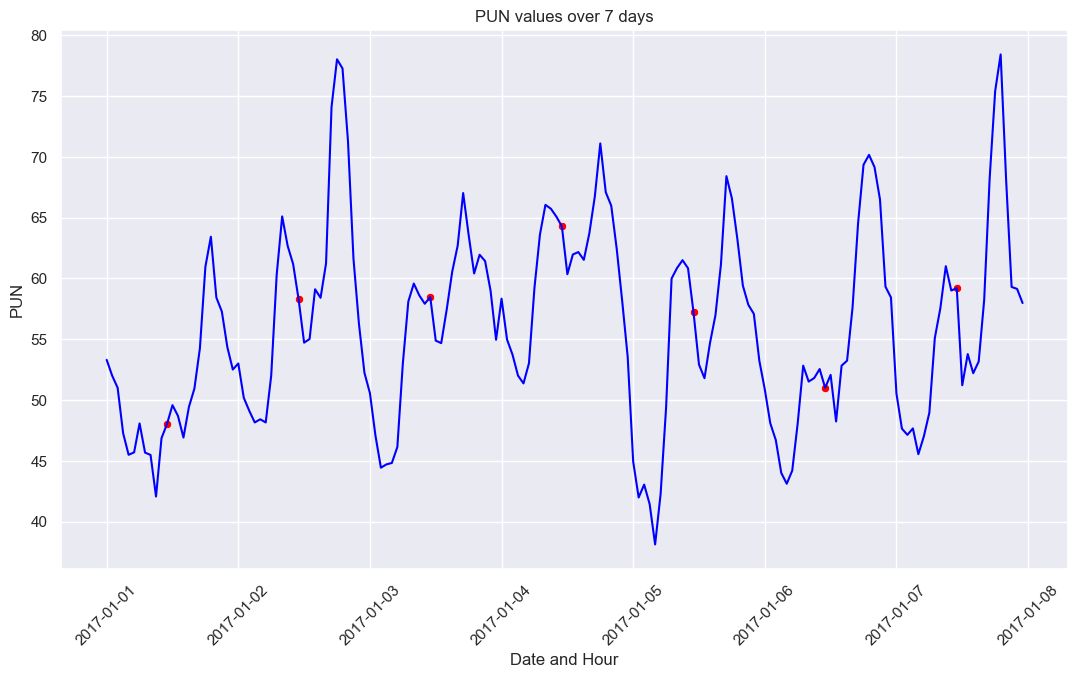}
    \caption{A batch of $7$ days electricity prices describes the state of the RL algorithm. The red dots correspond to the specific hour we want to build the optimal curves}.
    \label{fig:7days_pun}
\end{figure}

For each mode of production $k$ corresponding to $K$ different sources, we have:
    
- \textit{Unitary Production Costs}: The production cost \( C^k_t \) for each of the \( K \) production modes at time \( t \).

- \textit{Maximum Dispatched Volumes \( D^k_t \)}: The maximum volume that can be produced for each of the \( K \) production modes at time \( t \).

\textbf{Action}. The actions correspond to offering curves like the one depicted in Fig. \ref{fig:step}. The steps are described by a couple of quantities/prices $(V_i, P_i)$ the supplier intended to purchase. The number of steps $I$ the curves are constructed is a hyper-parameter of the model. Each {\em step} $i$ characterized by:
 
  - \textit{Volumes} \( V_i \): The volumes offered for sale.
  
  - \textit{Prices} \( P_i \): The prices at which these volumes are offered.

\begin{figure}[H]
    \centering
    \includegraphics[width=.91\textwidth]{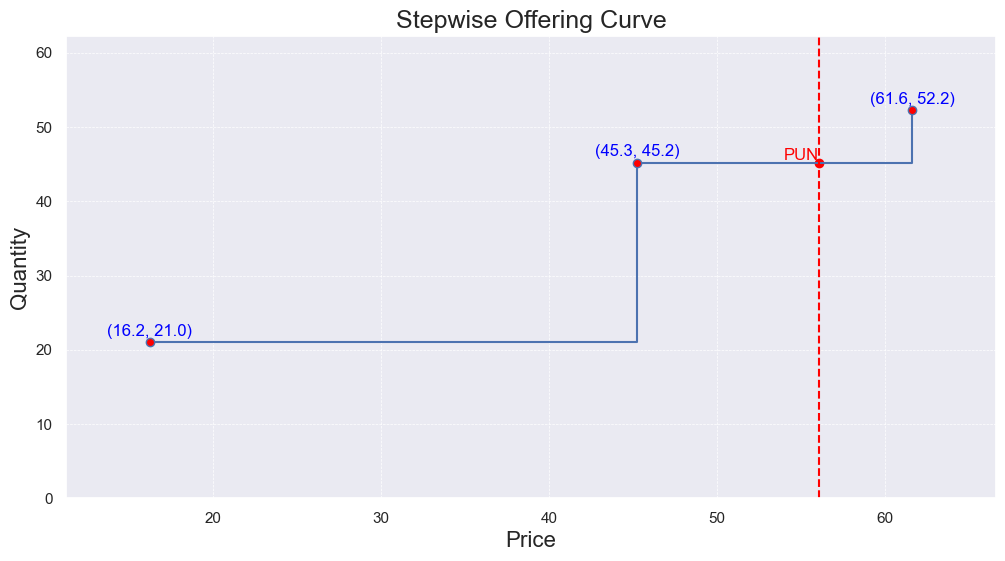}
    \caption{Offering curve corresponding to an output of the Actor NN. The reward is computed with Eq. \eqref{reward} using only the accepted offers, i.e. the ones with an offered price less than the registered pun (the red vertical line in the image).}
    \label{fig:output}
\end{figure}

  In Fig. \ref{fig:output}, we anticipate the space offering curve (we set $I=3$) obtained as outputs of the Actor-Network.

    The curve of the Actor NN is obtained as the offering that minimizes the reward (in terms of obtained profit) given the batch of prices for the last $7$ days.

\textbf{Reward}. The reward is computed by considering the volumes for which the offered price \( P_i \) is below the market clearing price \( PUN_{t+1} \). The revenue from the accepted offers (\( P_i \cdot V_i \)) is calculated, and then the production cost (\( C^i \cdot V_i \)) is subtracted, resulting in the final reward.

 Considering \( K \) modes of production and \( I \) bidding, i.e. steps for the offering curve, at each hour \( t \), the reward \( r_t \) is computed as:

\begin{equation}
    \label{reward}
    r_t (PUN_{t+1}, P_i, V_i, C_k, D_k) = \sum_{k=1}^K \sum_{i=1}^I \left( P_i \cdot V_i - C_k \cdot V_i \right) \mathds{1}_{\{ P_i \leq PUN_{t+1} \}}
\end{equation}
Where:

- \( PUN_{t+1} \) is the market clearing price for the same hour of the next day.

- \( P_i \) is the price offered at bidding \( i \).

- \( V_i \) is the volume offered at price \( P_i \).

- \( C_k \) is the unitary production cost for mode \( k \).

- \( D_k \) is the maximum producible volume for production mode \( k \).

- $ \mathds{1}_{\{ P_i \leq PUN_{t+1}\}} $ is the indicator function that equals 1 if the offered price \( P_i \) is less than or equal to \( PUN_{t+1} \), and 0 otherwise. 

\medskip

The presence of the indicator function ensures that the reward is computed by considering the revenue from the accepted offers (where \( P_i \leq PUN_{t+1} \)). According to Eq. \eqref{reward}, the reward corresponds to the realized profit,  computed by subtracting to the daily gain the production cost $C_k$ for the sold volumes $V_i$.

\medskip

We refer to Sec. \eqref{sec:ddpg} for a detailed study of how these quantities are merged into the DDPG algorithm, as well as other numerical simplifications introduced in the model.

\section{The DDPG Algorithm for Electricity Auction}
\label{sec:ddpg}

DDPG is a model-free actor-critic algorithm based on the deterministic policy gradient with continuous action spaces.

\subsection{The DDPG Algorithm}
 
In general terms, the DDPG algorithm consists of the following steps:

\begin{enumerate}
    \item Initialize the actor network with weights $\theta_\mu$ and the critic network with weights $\theta_Q$.

    \item Initialize the target networks for the actor and critic with the same weights: $\theta_\mu^{'} = \theta_\mu$ and $\theta_Q^{'} = \theta_Q$.
    
    \item Sample a minibatch of transitions $(s_t, a_t, r_t, s_{t+1})$ from the replay buffer.
    \item Update the critic network by minimizing the loss:
    \begin{equation}
    \label{criticloss}
    L = \frac{1}{N} \sum_{t} (y_t - Q(s_t, a_t | \theta_Q))^2
    \end{equation}
    where $y_t = r_t + \gamma Q(s_{t+1}, \mu(s_{t+1} | \theta_\mu^{'}) | \theta_Q^{'})$ is the target Q-value, and $\mu$ is the deterministic policy from the actor network. The critic loss \eqref{criticloss} corresponds to the Mean Squared Error (MSE) between the predicted Q-values and the actual rewards plus the discounted Q-value of the next state (i.e., the Bellman equation \eqref{bellman}). This loss guides the critic network to approximate the true Q-values \eqref{q-values} as accurately as possible.

    \item Update the actor-network using the sampled policy gradient introduced in Eq. \eqref{policy};
    \item Update the target networks using the soft update rule:
    \begin{equation}
    \begin{aligned}
    \theta_\mu^{'} \leftarrow \tau \theta_\mu + (1 - \tau) \theta_\mu^{'},  \quad
    \theta_Q^{'} \leftarrow \tau \theta_Q + (1 - \tau) \theta_Q^{'} \, ,
    \end{aligned}
    \end{equation}
    where $\tau \ll 1$ is a small constant that controls the update rate.

        \end{enumerate}

    An exploration strategy based on noise processes, e.g., the Ornstein-Uhlenbeck one, adds temporally correlated noise to the actions.
Accordingly, the noise process $X_t$ is defined by the following SDE:

\begin{equation}
\label{ou}
dX_t = -\theta (X_t - \mu) dt + \sigma dW_t \, , 
\end{equation}
that we discretize by the following Euler-Maruyama method:

\begin{equation}
\label{discr_ou}
X_{t+\Delta t} = X_t - \theta (X_t - \mu) \Delta t + \sigma \sqrt{\Delta t} \xi_t
\end{equation}

Where $\xi_t$ represents a random sample from a standard normal distribution.

Finally, the generated noise is added to the actions produced by the Actor-network:

\begin{equation}
\label{act_update}
a_t = \mu(s_t | \theta_\mu) + X_t
\end{equation}

where $a_t$ is the action taken at time $t$, $\mu(s_t | \theta_\mu)$ is the action produced by the actor-network for the state $s_t$, and $X_t$ is the noise generated by the Ornstein-Uhlenbeck process.

\subsection{Adaptation of the DDPG Algorithm for Electricity Auction}

We consider a single agent setting of an energy operator that interacts with the market (environment) in discrete time steps.

In summary, at each time step $t$, the agent:

\begin{enumerate}
    \item receives an observation $x_t$ consisting in of $24$-hour PUN array prices $P_t$ of the last $d$ days. 

    \item generates an action $a_t$, i.e. a stepwise curve modelling the offering curve corresponding to the output of the actor-network;
    \item observes a feedback scalar reward $r_t$ corresponding to the agent's payoff. 

        \end{enumerate}

\begin{figure}[H]
\includegraphics[scale=0.101]{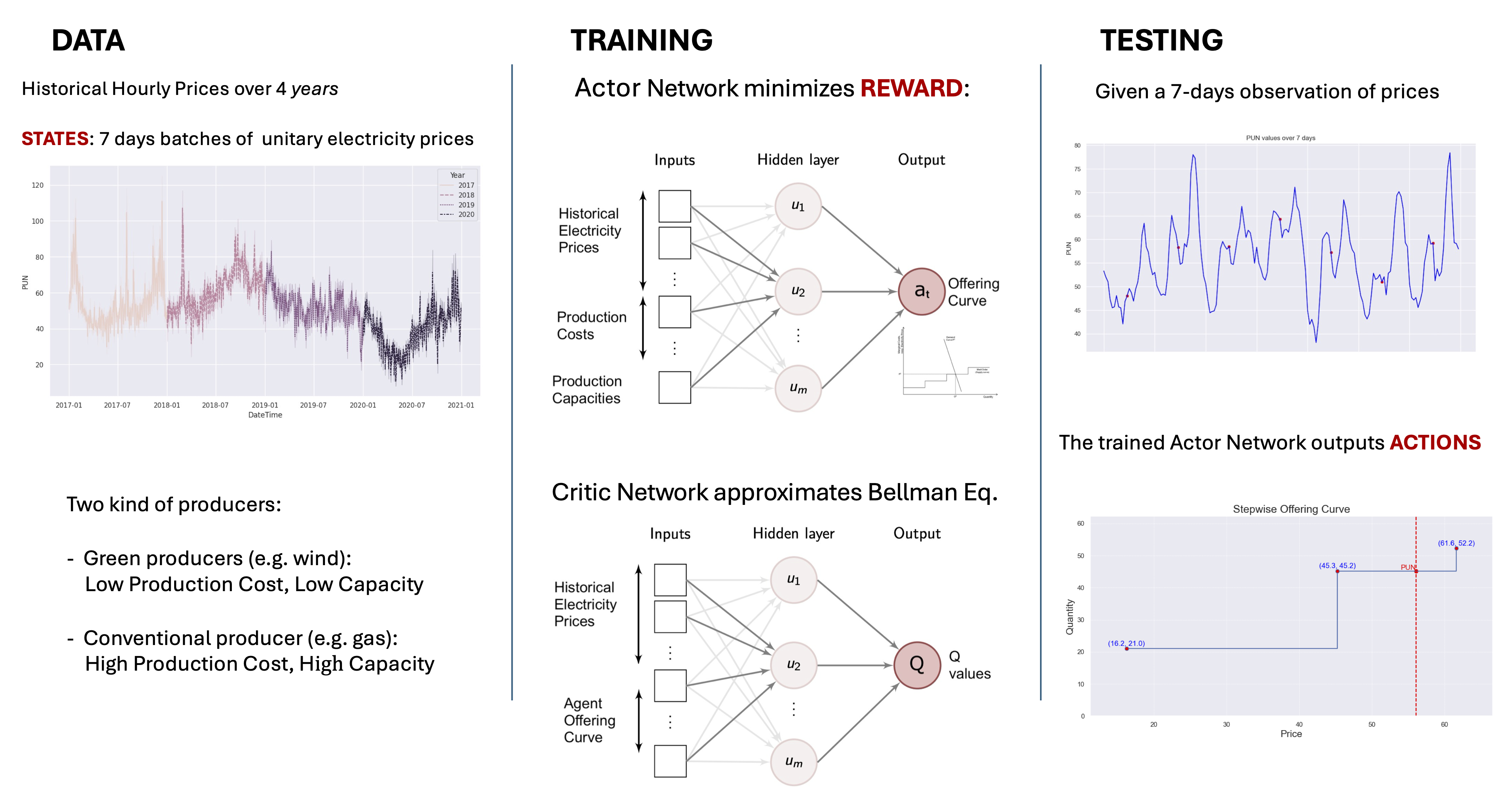}
    \caption{Graphical summary of the RL data-driven model.}
    \label{fig:arc}
\end{figure}

Both actor and critic are approximated using deep Feed Forward Neural Networks (FFNN) with a second set of target FFNNs.

\begin{figure}[H]
    \centering
    \includegraphics[width=.71\textwidth]{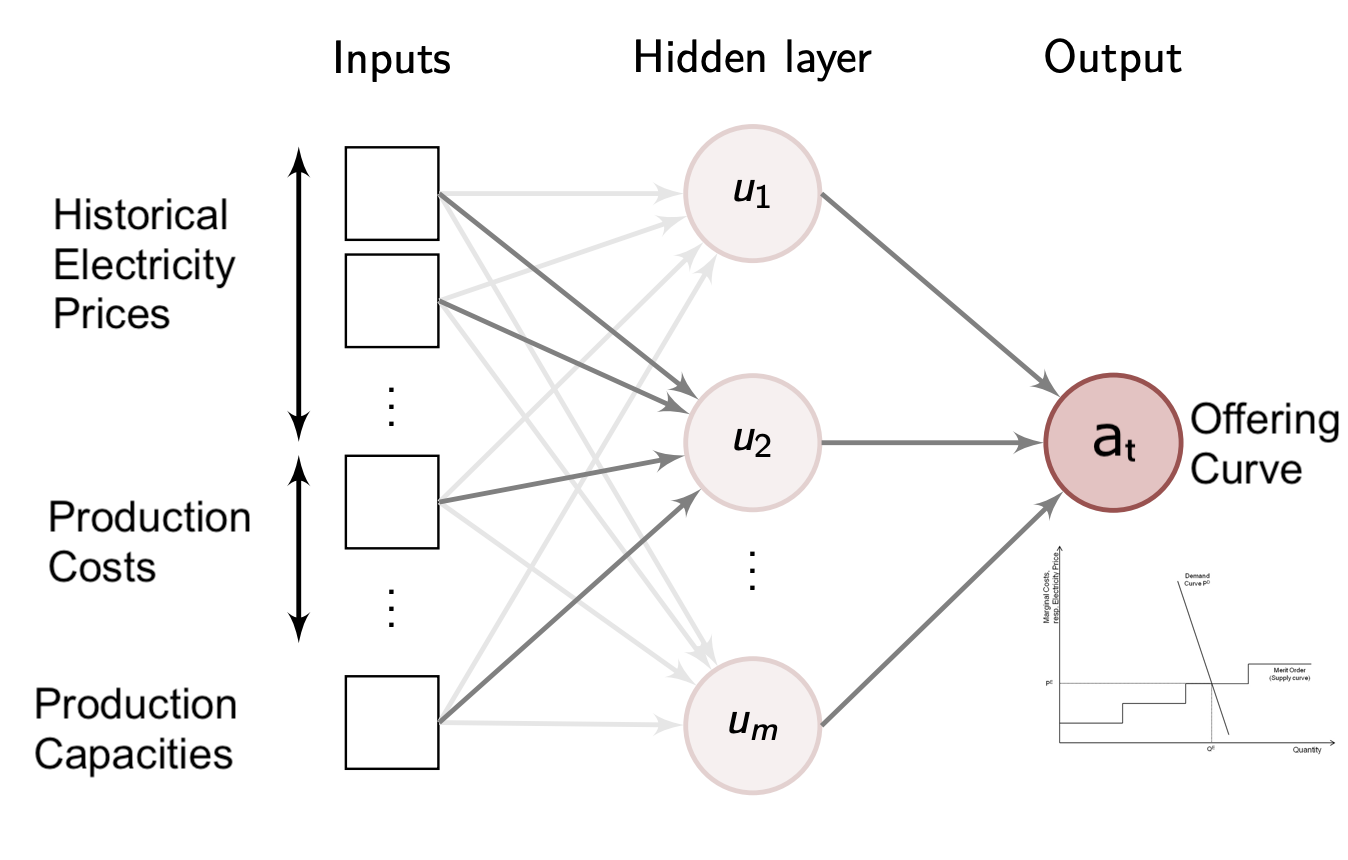}
    \caption{The Feed Forward Actor Network produces a vectorial output representing the Offering Curve.}
    \label{fig:actor}
\end{figure}

In Fig. \ref{fig:actor}, we sketch the Feed Forward NN we use to model the Actor Network.

The Critic Network approximates the Q function for a given state-action pair. On the other hand, the critic function $Q (s,a)$ approximates the Q-value value function given a pair (price-offering curve) that estimates the expected return by approximating the Bellman equation \eqref{bellman}.

\begin{figure}[H]
    \centering
    \includegraphics[width=.71\textwidth]{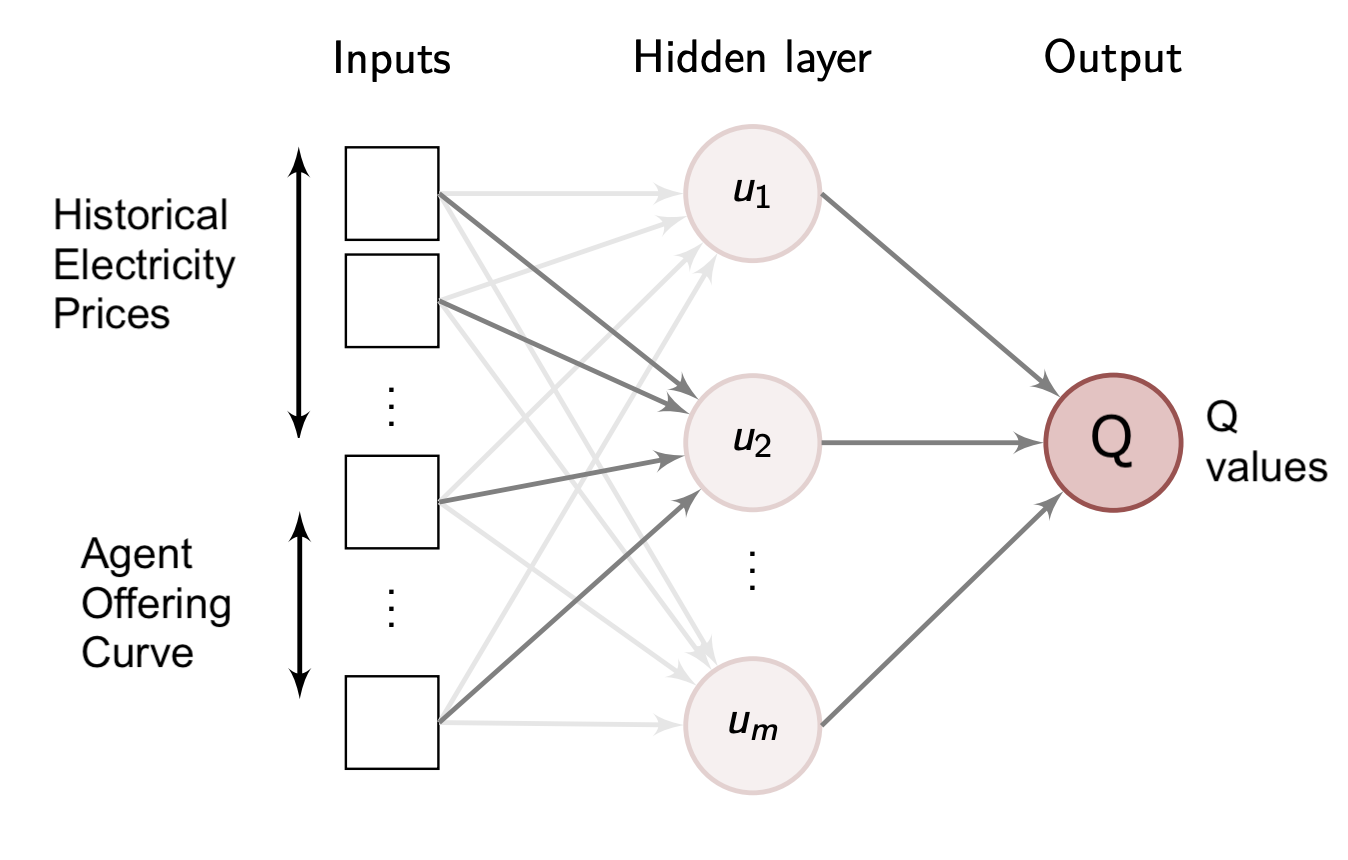}
    \caption{The Feed Forward Critic Network approximates the Bellman Equation given $(s_t, a_t)$}.
    \label{fig:critic}
\end{figure}

The actor is updated by applying the chain rule using the sampled policy gradient. 

\begin{equation}
\label{policy}
\nabla_{\theta_\mu} J \sim \frac{1}{N} \sum_t \nabla_a Q (s, a \vert \theta_Q) \big\vert_{\{s = s_t, a = \mu (s_t) \}}\nabla_{\theta_\mu} \mu (s \, \vert \, \theta_\mu) \big\vert_{\{s = s_t\}}
\end{equation}
that minimizes the distance between the current policy's actions and actions that maximize expected rewards. This technique estimates the gradient of the expected cumulative reward concerning the policy parameters using samples collected during interactions with the environment.

\section{Implementation of the algorithm}
\label{sec:numerics}

Concerning the quantities introduced in Sec. \ref{sec:model2}, we assume curves composed of $I=3$  biddings. We set $K=3$ modes of production and $I=3$ fixed generation capacities and different production costs corresponding to different energy sources. We consider a constant array of production costs and available power to calibrate our result to calculate a value for the reward function defined in Eq. \eqref{reward}. In particular, we set the number of sources of production $K$ to $3$ and set the cost $C^i = [10,30,60]$ and the production capacity $D^i = [30,200,800]$. This choice, from a modelling point of view, the source with a low marginal cost 
and low capacity, i.e. $(10, 30)$, represents a renewable source of production, the one with a high marginal cost and capacity a conventional one (e.g. gas) plus an intermediate one. Without loss of generality, one could also consider stochastic quantities, with the additional effort of storing the ad-hoc Stochastic Differential Equation (SDE) simulation while adding an extra source of randomness to the NN input. We decide to leave this additional feature for future work.


Another hypothesis we introduce is that the agent only knows about its expenses, available resources, historical electricity prices and nothing about its competitors. As a result, its bidding strategy can be represented as a stochastic process that adheres to a decision-making framework.

\medskip

Concerning the development algorithm, we report some tools we used in the DDPG method to improve the learning process's stability:

\begin{itemize}
    \item a replay buffer is employed to store past experiences and sample mini-batches of transitions (the so-called experiences arrays) for training;

    \item using different learning rates for actor and critic networks. The target networks are updated slowly, using a soft update rule with a small mixing factor, which helps to stabilize learning;

    \item Adding regularisation techniques, such as L2 regularisation, to the loss functions for the actor and critic networks helps prevent overfitting;

    \item carefully tuning the hyperparameters, such as the learning rates, discount factor, and soft update rate, significantly impacts the performance of DDPG. We conduct a systematic search or optimization techniques to find the best hyperparameters.

\end{itemize}

Besides working in mini-batches, we add common noise to randomize the actions. We assume noise as a discretized Ornstein-Uhlenbeck process as defined in Eq. \eqref{discr_ou}. By an empirical calibration, we set the value of the rate of mean reversion $\theta$ to $0.15$, the mean $\mu$ to $1$ with diffusion $\sigma$ taking values in the interval [$1$, $10$] guaranteeing an efficient impact for the update of the action \eqref{act_update}.

Consequently, we derive the following scheme:

\begin{algorithm}[H]
\caption{DDPG for electricity auctions}\label{alg:cap}
\begin{algorithmic}[1]
\State Initialize the Actor network $\mu(s_t | \theta_\mu)$ and the Critic network $Q(s_t, a_t | \theta_Q)$ with random weights $\theta_\mu$ and $\theta_Q$.
\State Initialize the target networks $\mu'(s_t | \theta_\mu')$ and $Q'(s_t, a_t | \theta_Q')$ with weights $\theta_\mu' \leftarrow \theta_\mu$ and $\theta_Q' \leftarrow \theta_Q$.
\State Initialize the Ornstein-Uhlenbeck noise process $X_t$.
\State For each episode:
\begin{algorithmic}[A]
\State Initialize the environment and obtain the initial state $s_0$.
\State For each time step $t$:
\begin{algorithmic}[1]
\State Select the action $a_t = \mu(s_t | \theta_\mu) + X_t$, where $X_t$ is the noise generated by the Ornstein-Uhlenbeck process.
\State Execute the action $a_t$ in the environment and observe the reward $r_t$ and the next state $s_{t+1}$.
\State Store the transition $(s_t, a_t, r_t, s_{t+1})$ in the replay buffer.
\State Update the noise process $X_{t+\Delta t} = X_t - \theta (X_t - \mu) \Delta t + \sigma \sqrt{\Delta t} \xi_t$, where $\xi_t$ is a random sample from a standard normal distribution.
\State If the replay buffer contains enough samples, sample a mini-batch of transitions $(s_j, a_j, r_j, s_{j+1})$ from the replay buffer.
\State Update the Critic network by minimizing the loss:
\begin{equation}
L(\theta_Q) = \frac{1}{m} \sum_{j=1}^m \left( Q(s_j, a_j | \theta_Q) - (r_j + \gamma Q'(s_{j+1}, \mu'(s_{j+1} | \theta_\mu') | \theta_Q')) \right)^2
\end{equation}
\State Update the Actor Network using the sampled policy gradient:
\begin{equation}
\nabla_{\theta_\mu} J(\theta_\mu) \approx \frac{1}{m} \sum_{j=1}^m \nabla_{a} Q(s_j, a | \theta_Q) \big|_{a=\mu(s_j | \theta_\mu)} \nabla_{\theta_\mu} \mu(s_j | \theta_\mu)
\end{equation}
\State Update the target networks using soft updates:
\begin{equation}
\begin{aligned}
\theta_\mu' \leftarrow \tau \theta_\mu + (1 - \tau) \theta_\mu', \quad
\theta_Q' \leftarrow \tau \theta_Q + (1 - \tau) \theta_Q',
\end{aligned}
\end{equation}
\end{algorithmic}
\end{algorithmic}
\State Repeat until the desired level of performance is achieved or a maximum number of episodes is reached.
\end{algorithmic}
\end{algorithm}

\subsection{Selection of the DDPG hyperparameters}

We briefly describe the hyperparameters we use for the DDPG algorithm while referring to Table \ref{tab16} for the reference values we set for the implementation.

\begin{table}[H]
\begin{center}
\caption{DDPG Hyperparamters}\label{tab16}%
\begin{tabular}{lc}
\toprule
{\em Episodes}   & $1000$, $1500$  \\
\midrule  {\em Length of an Episode}  [Days]     & $15$, $20$, $30$  \\
\midrule {\em Batch Size}    & $64$ \\
\midrule {\em Hidden Size} \hspace{4cm}  & $64$    \\
\midrule {\em Actor Learning Rate}  & $0.0001$   \\
\midrule {\em Critic Learning Rate}    & $0.00001$       \\
\midrule {\em Discount Factor}    & $0.99$  \\
\midrule {\em Tau}  & $0.01$ \\
 \midrule {\em Max Memory Size}     & 50000  \\
\bottomrule
\end{tabular}
\end{center}
\label{hyper}
\end{table}

An episode means dealing with a complete sequence of interactions between the agent and the environment. The {\em Episodes} corresponds to the times the agent will engage with the environment to learn and improve its policy. The {\em Length of an Episode} refers to the number of time steps or interactions the agent experiences within a single episode. It defines how long the agent operates in the environment before the episode concludes. The {\em Batch Size} refers to the number of experiences (state action-reward-next state tuples) sampled from the replay buffer at each iteration of the training process. The {\em Hidden Size} refers to the number of neurons or units in the hidden layers of the neural networks used in the DDPG algorithm. The {\em Actor Learning Rate} and the {\em Critic Learning Rate} control step size of the gradient update at which the actor network's weights $\theta_\mu$ and the critic network's weights $\theta_Q$ are updated during training. The {\em Discount Factor} parameter $\gamma$ introduced in Eq. \eqref{discount} represents the relative importance of future rewards compared to immediate rewards. A higher gamma value places more importance on long-term rewards, potentially encouraging the agent to consider the future consequences of its actions. The {\em Target Update Factor} $\tau$ is a relaxation factor that defines how often the parameters are copied from the original networks to the target for the network parameters copied. The {\em Max Memory Size} determines the capacity of the replay buffer, which is a crucial component in DDPG. The replay buffer stores past experiences (state, action, reward, next state) that the agent uses to learn from.

Since DDPG hyperparameters are interconnected, finding the correct tuning significantly impacts the algorithm's performance and stability. We calibrate the values of Table \ref{tab16} by an experimental procedure by adjusting them to optimize the global reward. We refer to the following Subsection for a complete overview of the tuning procedure and some technical considerations.




\subsection{Numerical Results}

\begin{figure}[H]
    \centering
    \includegraphics[width=.85\textwidth]{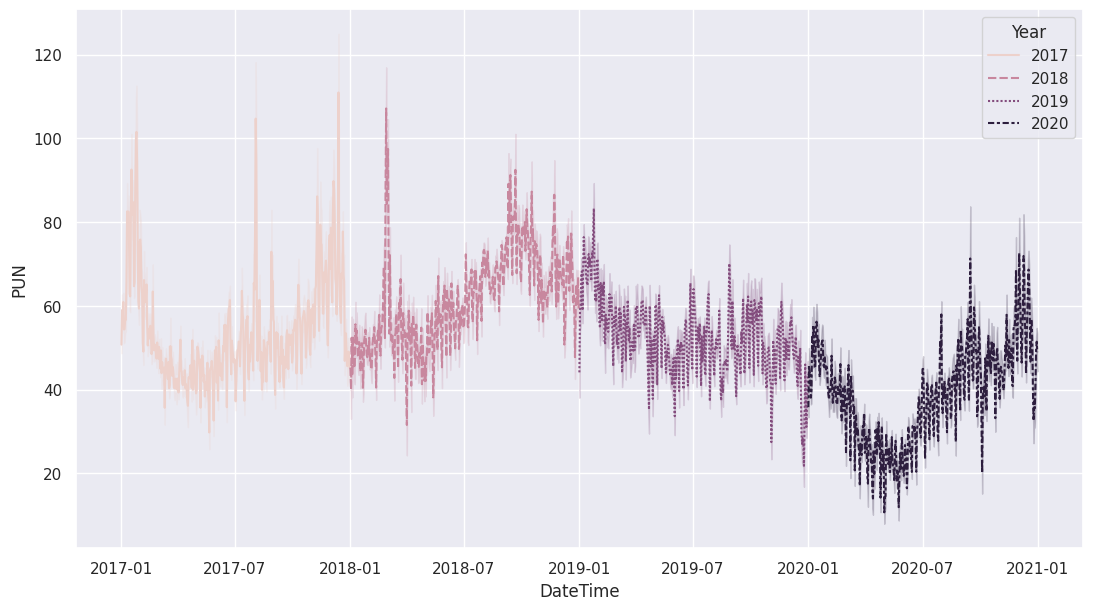}
    \caption{Italian hourly PUN from 01-Jan-2017 to 12-Dec-2020 for a total of $35064$ data points }.
    \label{fig:pun}
\end{figure}

The dataset consists of Italian electricity prices over 4 years. Precisely, we use hourly data for PUN from 2017 to 2020 that we download at \url{https://www.mercatoelettrico.org/it/Download/DatiStorici.aspx} plotted in Fig. \ref{fig:pun}.

We use stratified sampling to ensure that training and testing sets have a representative mix of data throughout the timeframe.

Rewards obtained from the environment might vary significantly since they depend on the price time series that is highly not stationary, as we can see from Fig. \ref{fig:pun}. To address this issue, we include a normalization factor into the reward to consider the potential {\em normalized reward} that can be obtained each time $t$. 

The normalized reward, denoted as $\mathcal{R}_{norm}$, is calculated by dividing the actual reward $r$ introduced in Eq. \eqref{reward} for a given time step by the maximum possible value of the reward $r_{max}$ that can be achieved under the particular time step conditions. When production costs, capacity, and the matching pun for a specific time step are considered, the maximum reward reflects the most profit that may be realized. For more steady and practical learning, this normalization scales the reward values to a consistent range between 0 and 1: when the agent receives the maximum reward, the normalized reward will be 1. The normalized reward ranges from 0 to 1 if the agent performs below the maximum. This option enables the agent to focus on tailoring its strategy in response to the relative performance improvement, which improves the consistency of comparisons between various learning contexts.

\medskip

Besides the normalized reward, we plot two interesting metrics: policy loss and critical loss. The agent learns to improve its policy by adjusting the weights of the actor-network to maximize the expected cumulative reward $J$, defined in the equation. \eqref{functinoal}. The policy loss measures the discrepancy between the actions chosen by the current policy and those that would lead to higher expected rewards. It is approximated by the sampled policy gradient introduced in \eqref{policy}. The Critic Loss introduced in Eq. \eqref{criticloss}
is associated with training the critic network, measures the precision of the predictions of the Q-value of the critic and guides the critic network to approximate the Q-values that satisfy the Bellman equation \eqref{bellman}.

We point out that the algorithm is sensitive to the initialization of parameters such as ($\#$ episode, $\#$ days in the episode, production cost) that must be carefully chosen.

For this reason, we perform several simulations for different values of {\em Episodes} and {\em Length of an Episodes}, reporting the obtained results for 5 different simulations for which we get the following plots.

\begin{figure}[H]
\centering
\begin{tabular}{cc}
\multicolumn{2}{c}{\includegraphics[width=.7\textwidth]{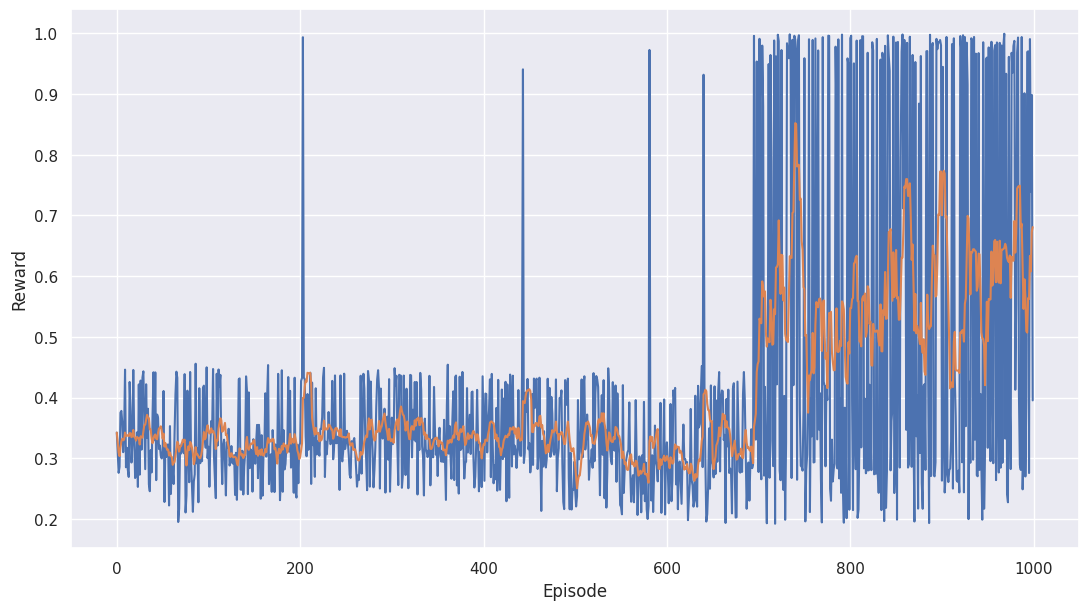}} \\
\includegraphics[width=.4\textwidth]{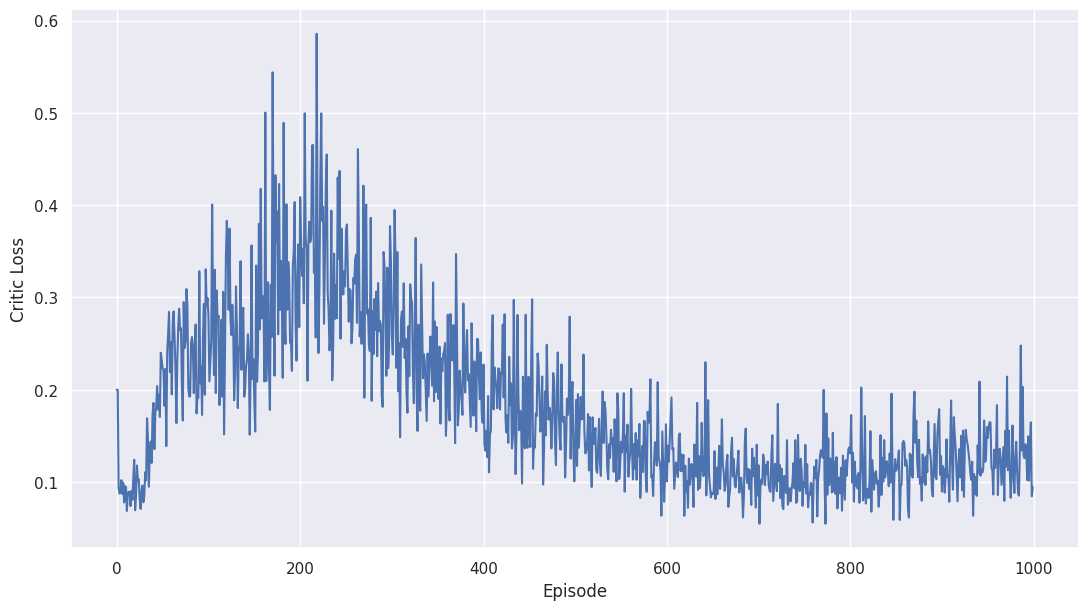}&
\includegraphics[width=.4\textwidth]{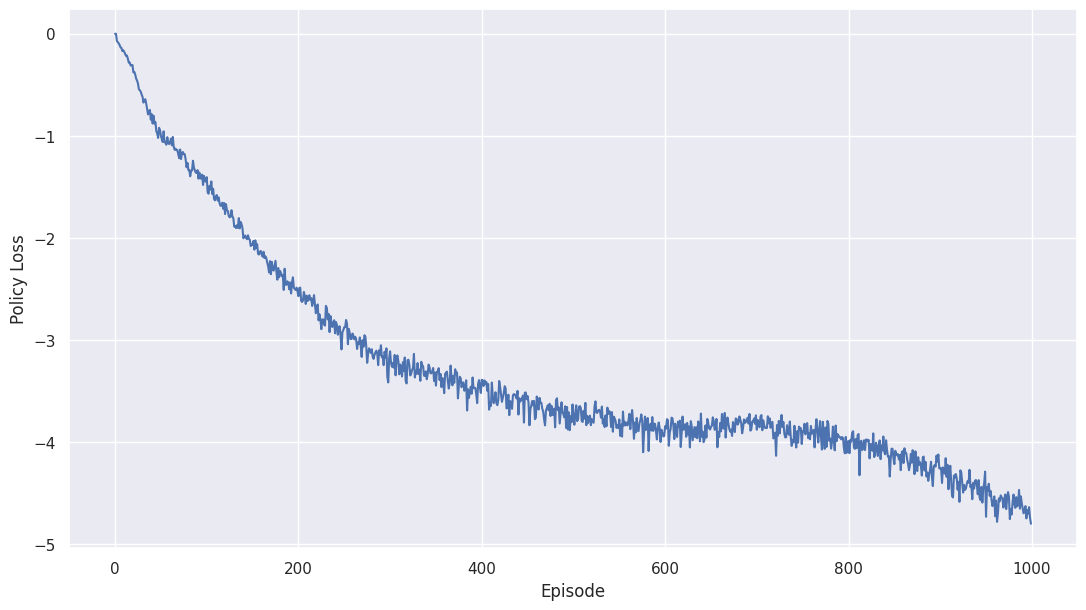}\\
\end{tabular}
    \caption{{\bf Simulation 1}: Episodes = 1000, length of episode = 30}
    \label{validation1} 
\end{figure}

\begin{figure}[H]
\centering
\begin{tabular}{cc}
\multicolumn{2}{c}{\includegraphics[width=.66\textwidth]{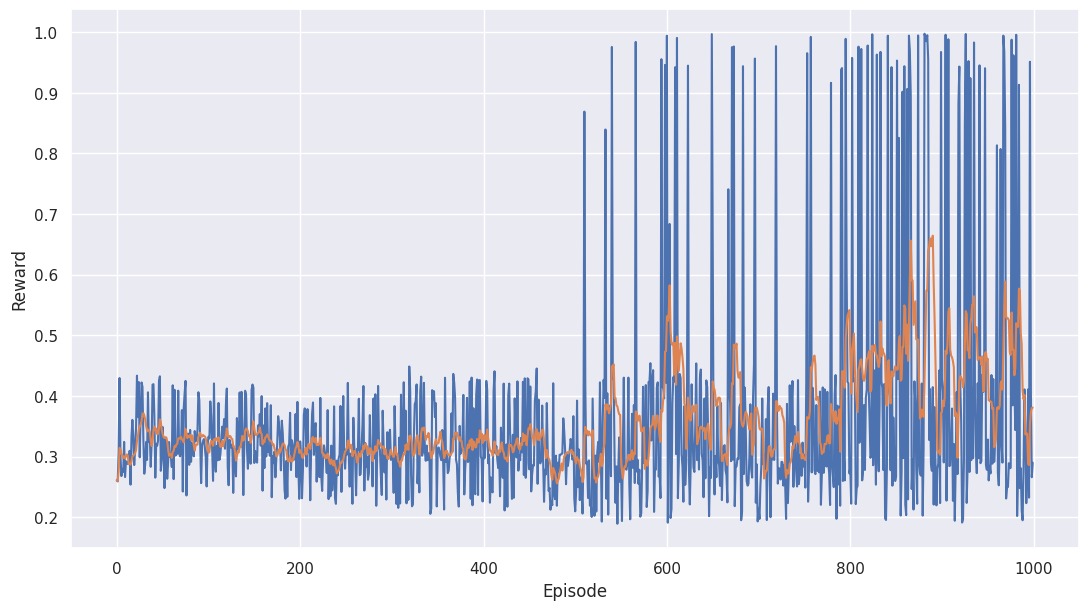}} \\
\includegraphics[width=.38\textwidth]{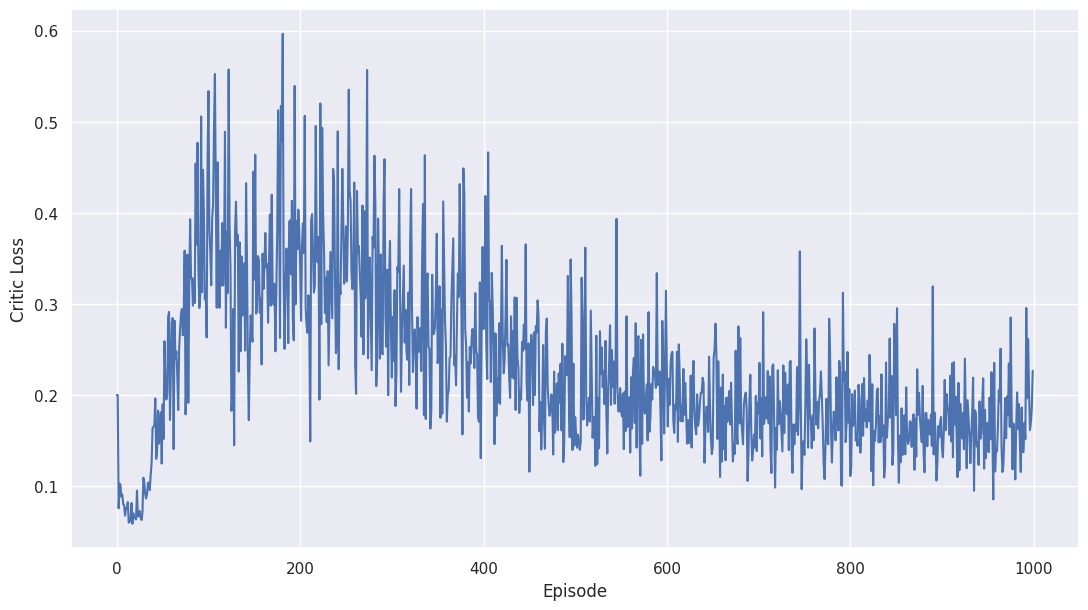}&
\includegraphics[width=.38\textwidth]{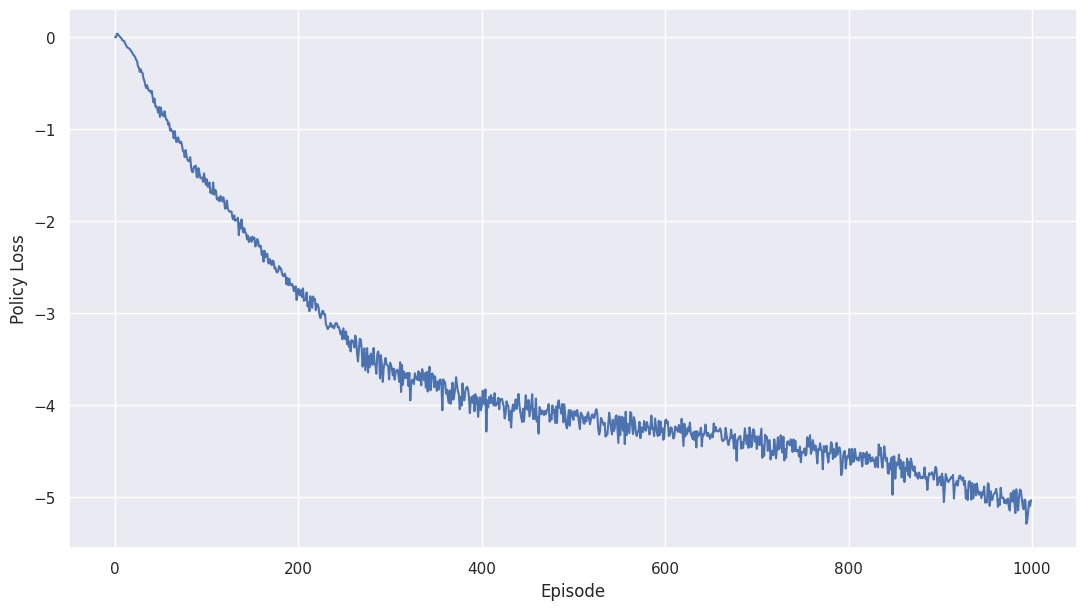}\\
\end{tabular}
    \caption{{\bf Simulation 2}: Episodes = 1000, length of episode = 30}
    \label{validation2} 
\end{figure}

\begin{figure}[H]
\centering
\begin{tabular}{cc}
\multicolumn{2}{c}{\includegraphics[width=.66\textwidth]{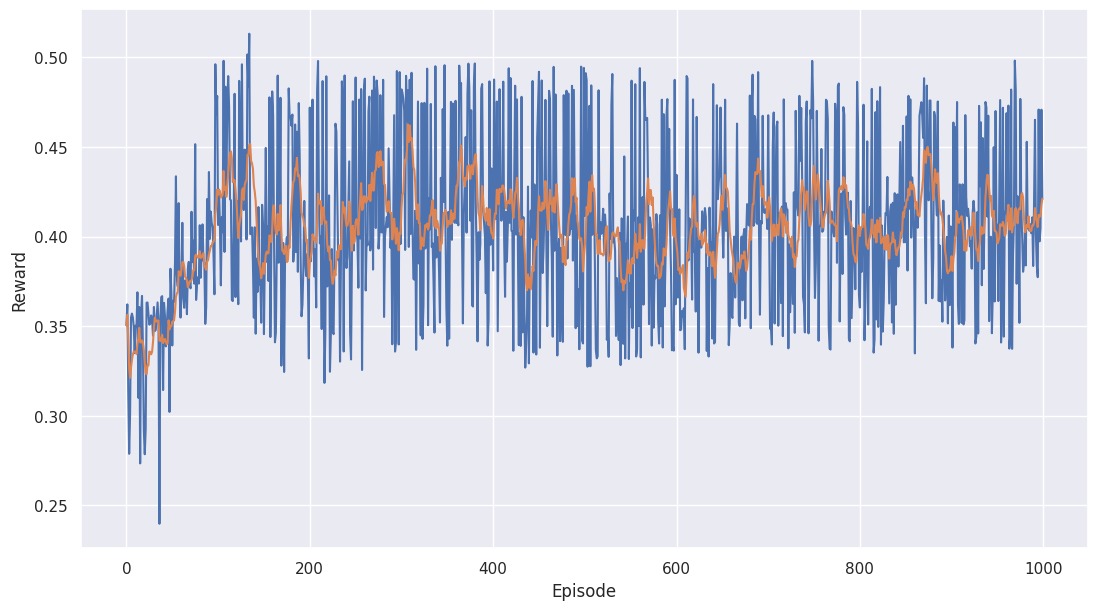}} \\
\includegraphics[width=.38\textwidth]{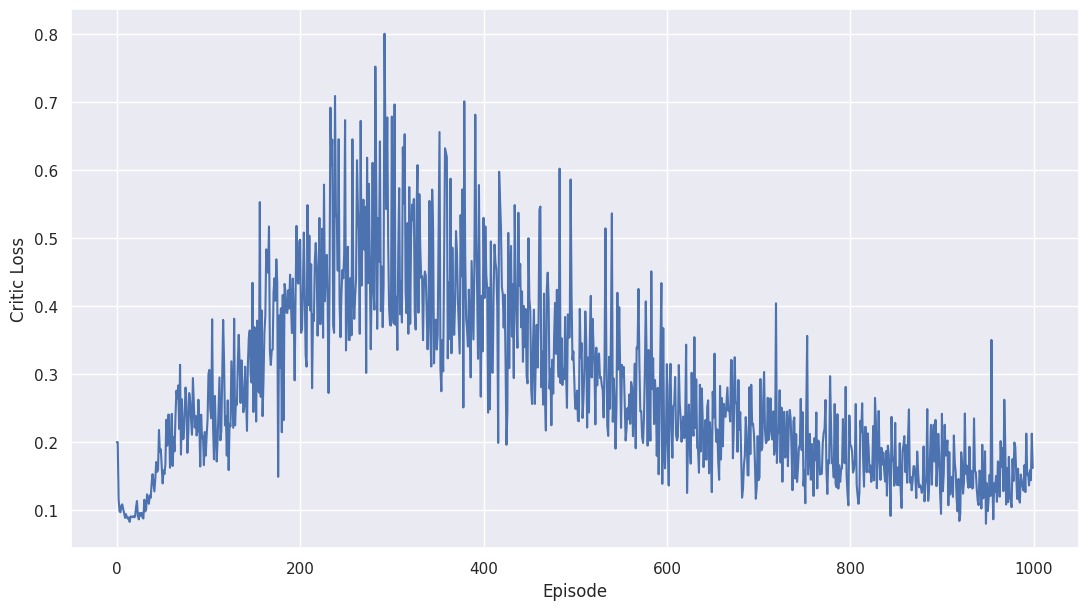}&
\includegraphics[width=.38\textwidth]{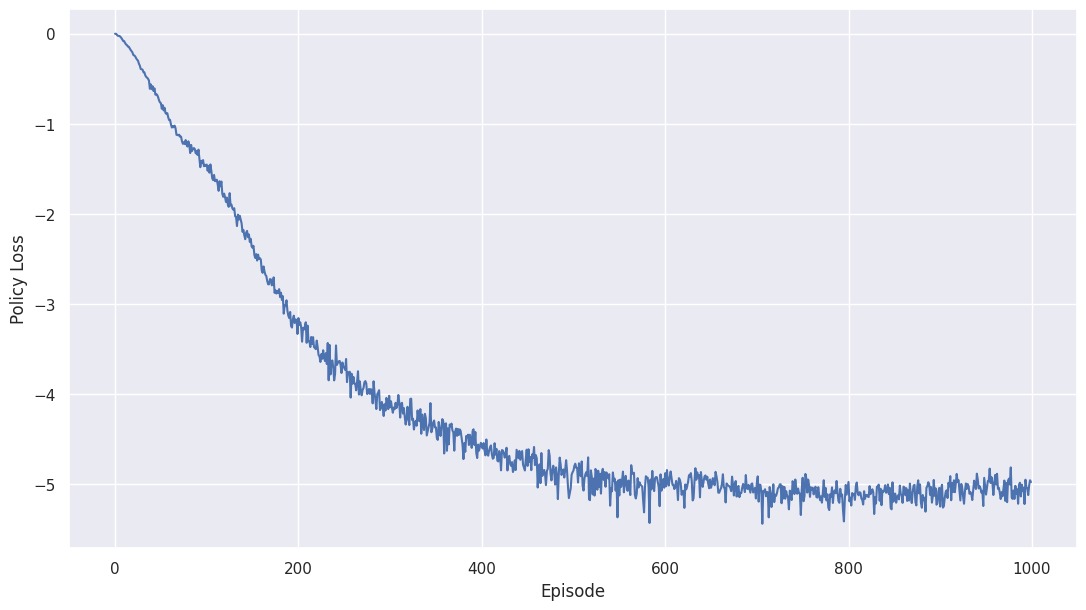}\\
\end{tabular}
    \caption{{\bf Simulation 3}: Episodes = 1000, length of episode = 30}
    \label{validation3} 
\end{figure}

\begin{figure}[H]
\centering
\begin{tabular}{cc}
\multicolumn{2}{c}{\includegraphics[width=.66\textwidth]{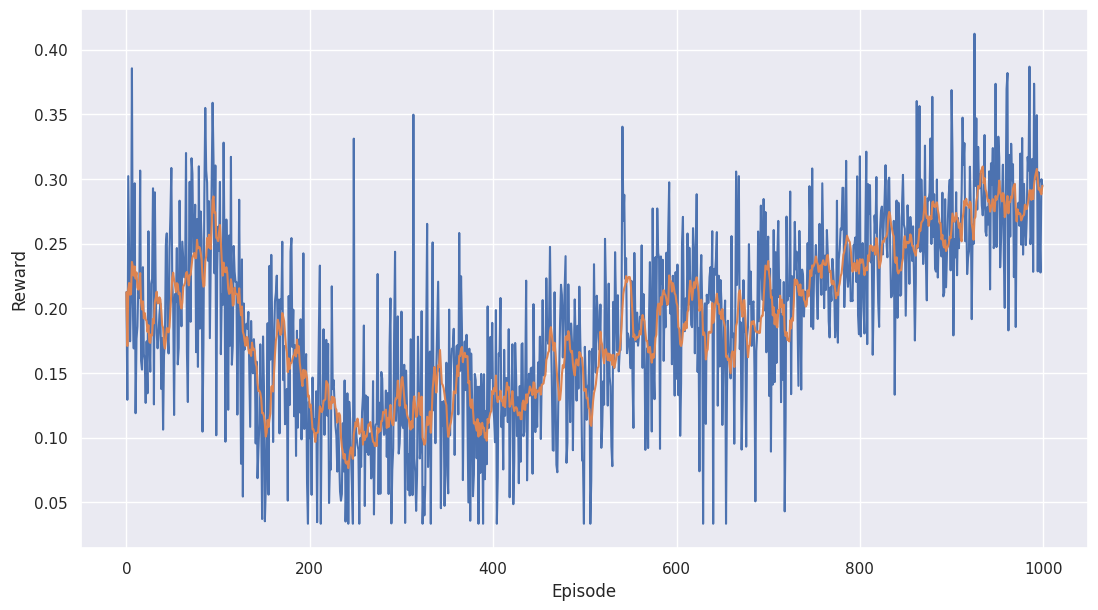}} \\
\includegraphics[width=.38\textwidth]{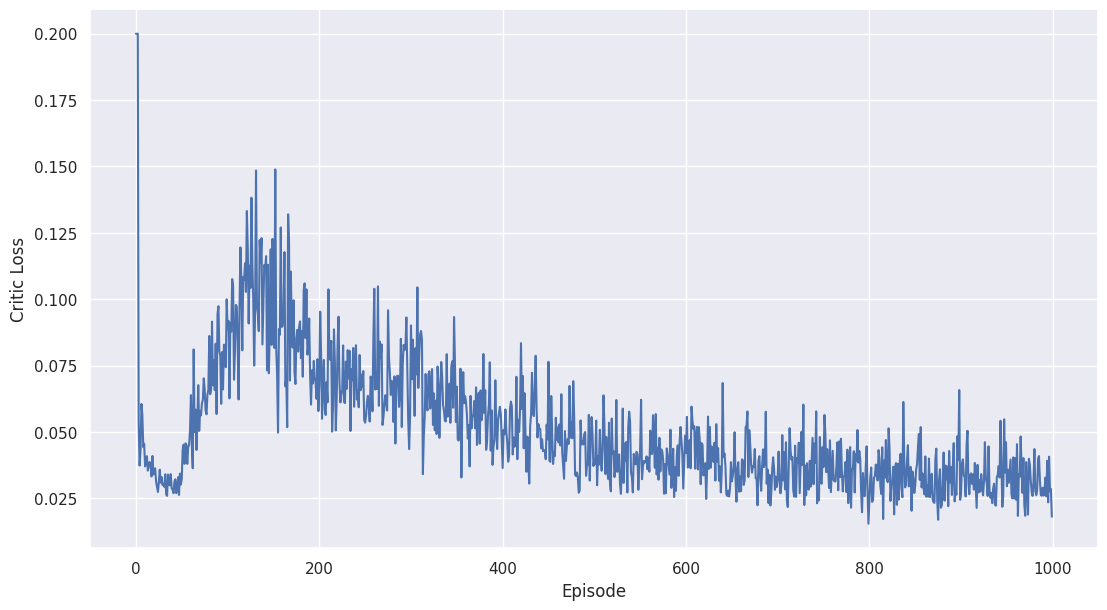}&
\includegraphics[width=.38\textwidth]{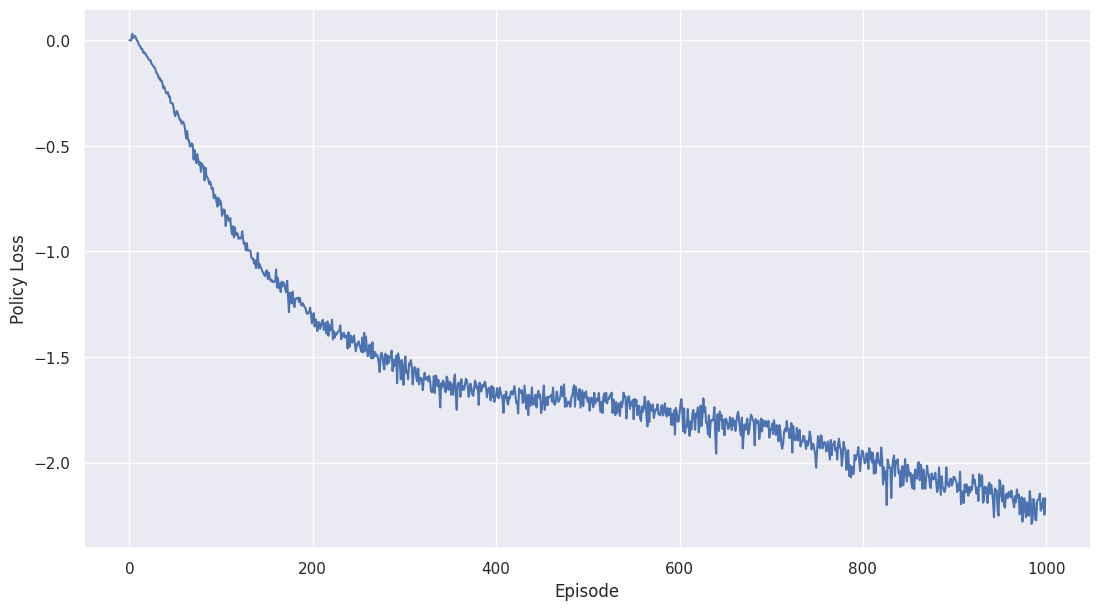}\\
\end{tabular}
    \caption{{\bf Simulation 4}:Episodes = 1000, length of episode = 20}
    \label{validation4} 
\end{figure}

In all plots, we can see how the policy loss consistently decreases over time; hence, the actor moves towards a better policy. In contrast, the critic loss shows an initial spike before gradually decreasing. Thus, the Q-value estimate improves only after approximately $300$ {\em Episodes.} This behaviour is a typical pattern in the training dynamics of RL algorithms that can be attributed to the interaction with a stochastic environment. Specifically, in the early stages of training, the agent's policy might be far from optimal, leading to higher Q-values and a more considerable critic loss. As the agent explores the environment and gathers more experience, it gradually refines its policy, causing the loss of the critic to decrease. Moreover, this delayed learning is also linked to the target networks, which are updated slowly, leading to a more significant initial critical loss. The critic loss decreases as training progresses and the target networks catch up.

However, the gradual decrease in policy and critic losses as training metrics and the increase of the normalized reward indicate the agent's convergence towards more optimal policies.

\subsection{Considerations, Challenges and Limitations of the model}

As previously encountered in other settings such as \cite{markov}, the market rewards may be highly stochastic in the single-agent framework. They may be unable to converge to a fixed point, leaving to the oscillatory behaviours of market producers. This behaviour is expected because generators do not account for the strategies of their competitors, which may evolve. In contrast, the gradual decline in policy loss and the initial spike and subsequent decrease in critic loss are signs of the DDPG algorithm's learning process. They demonstrate the agent's evolution from initial exploration and poorly predicted Q-values to more precise value estimations and a more developed strategy. These actions indicate that the algorithm is learning from its training data and adjusting to its surroundings, which are generally expected behaviours.

\bigskip

The difficulties we run across while building the algorithm can be divided into two groups: those caused by how the electrical auction markets are set up and those specifically related to the RL algorithm. The first type consists of the action space's partial observability (the agent can learn about other people's actions and tactics only through time series of historical pricing), non-stationarity of the data, and high problem dimensionality.
 Other challenges sources arising from the stability of the DDPG are those connected to inherent difficulties in learning, such as function approximation errors and non-stationary targets.

Additionally, when we monitor the algorithm's process, we notice the following observations, which are primarily concerned with the stability of the DDPG technique and the selection of the hyperparameters:

\begin{itemize}
    \item A higher number of episodes allows the agent to explore the environment more extensively, potentially leading to better policy convergence. Episodes provide opportunities for the agent to encounter a diverse range of states and price profiles, which is crucial for learning a robust policy that performs well across different scenarios. Conversely, the learning process often exhibits diminishing returns with increasing episodes. Initially, the agent might rapidly learn and improve its policy, but over time, the rate of improvement might slow down as it explores less novel situations. Summarizing, the length of an episode appears to be a crucial parameter for the exploration-exploitation trade-off;

\item  Also {\em Length of an Episode} hyperparameter can influence the agent's learning dynamics and exploration strategy. A more extended episode provides the agent more time to investigate the surroundings thoroughly and make more choices, which might result in further exploration. A shorter episode, on the other hand, would encourage the agent to focus more on using what it already knows, thereby limiting exploration. We identify a crucial value of $15$ days for {\em Length of an Episode} below which we do not have convergence;

\item   By experimentation, we set the learning rate of the critic lower than the one of the actor-network. Updates to obtain more stable training by preventing one network from significantly outpacing the other;

\item A larger batch size can lead to more efficient learning as it allows the agent to learn from more experiences in parallel. However, we note that employing a huge batch size may increase the amount of noise and volatility in the learning process, resulting in less consistent training. As each update is based on a smaller subset of the overall experiences, we selected a batch size that is relatively modest to encourage the agent to explore more varied encounters.

\end{itemize}

\section{Conclusions and possible future directions}

\label{sec:future_directions}

The article's goal concerns the development of an optimizing strategy for a single-agent RL setting. We develop a DDPG algorithm to learn a deterministic policy through actor-critic architecture using a Q-value function to grade the proposed offering curve. This combination allows DDPG to effectively handle continuous action spaces, making it a robust algorithm for various applications.

We extend the result in \cite{xho} by introducing NNs, precisely the DDPG method, and by directly using historically accurate prices during decision-making. 

We conclude the article by presenting some extensions of this model that may be worth further investigating:  we solely test feed-forward NNs for implementing this model. We think the algorithm's performance could be enhanced by choosing architectures that are better suited for time series, including recurrent NNs, such as Long Short Term Memory NNs;
 we model production cost and production capacity as deterministic constant quantities. Incorporating stochastic capacities to model random fluctuations associated with renewable energy production could provide significant benefits. This research serves as a foundational building block for developing algorithms for multi-agent or Mean Field systems, for example, by integrating the empirical distribution of the offering curves of other operators. Thereafter, the construction of a distributed optimization system -- an algorithm based on decentralized coordination such as local rewards or consensus schemes -- is the foundation for extending a single agent into a multi-agent context.


 \bibliographystyle{elsarticle-num} 



\end{document}